\documentclass[11pt]{amsart}
\usepackage{amsmath}
\usepackage{amssymb}
\usepackage{amsthm}
\usepackage{latexsym}
\usepackage{graphicx}
\usepackage{hyperref}
\usepackage{enumerate}

\newtheorem{thm}{Theorem}[section]

\newtheorem{lem}[thm]{Lemma}

\newtheorem{rem}[thm]{Remark}
\newtheorem{thmintro}{Theorem}

\providecommand{\norm}[1]{\left\| #1 \right\|}
\newcommand{\enuma}[1]{\begin{enumerate}[\textup{(}a\textup{)}] {#1} \end{enumerate}}

\newcommand{\mh}{\mathbb}
\newcommand{\mr}{\mathrm}
\newcommand{\mc}{\mathcal}
\newcommand{\mf}{\mathfrak}

\newcommand{\ts}{\textstyle}

\newcommand{\Z}{\mathbb Z}

\newcommand{\R}{\mathbb R}
\newcommand{\C}{\mathbb C}
\newcommand{\ep}{\epsilon}

\newcommand{\af}{\mr{aff}}

\newcommand{\inp}[2]{\langle #1 \,,\, #2 \rangle}

\begin{document}

\title{A formula of Arthur and affine Hecke algebras}

\author{Eric Opdam}
\address{Korteweg-de Vries Institute for Mathematics\\
Universiteit van Amsterdam\\
Science Park 904\\
1098 XH Amsterdam\\
The Netherlands}
\email{e.m.opdam@uva.nl}
\author{Maarten Solleveld}
\address{Mathematisches Institut\\
Georg-August-Universit\"at G\"ottingen\\
Bunsenstra\ss e 3--5\\
37073 G\"ottingen\\
Germany}
\email{maarten@uni-math.gwdg.de}
\date{\today}

\maketitle

\begin{abstract}
Let $\pi, \pi'$ be tempered representations of an affine Hecke algebra
with positive parameters. We study their Euler--Poincar\'e pairing $EP (\pi,\pi')$,
the alternating sum of the dimensions of the Ext-groups. We show that $EP (\pi,\pi')$
can be expressed in a simple formula involving an analytic R-group, analogous to
a formula of Arthur in the setting of reductive $p$-adic groups. Our proof applies
equally well to affine Hecke algebras and to reductive groups over
nonarchimedean local fields of arbitrary characteristic.
\end{abstract}

\tableofcontents

\section*{Introduction}
Let $\mathbb F$ be a non-archimedean local field, and
let $L$ be the group of rational points of a connected reductive
algebraic group $\mathbf{L}$ defined over $\mathbb F$.
Affine Hecke algebras are powerful tools to describe
Bernstein blocks \cite{BeDe} of the category of smooth representations
of $L$ explicitly, in a way compatible with harmonic analysis.
For example, if a Bernstein block $\mc{B}$ admits a type $\tau$ then the
theory of types \cite{BK} establishes a Morita
equivalence between $\mathcal{B}$ and the module category
of the Hecke algebra $\mathcal{H}_\tau$ of the type $\tau$, such that the
Plancherel measure of $L$ restricted to $\mathcal{B}$ corresponds to the
spectral measure of $\mathcal{H}_\tau$ viewed as a Hilbert algebra \cite{BHK}.
In various cases the Hecke algebras of types have been described
explicitly as affine Hecke algebras \cite{Mor,Lus-Unip}. Other
categorical equivalences describing Bernstein blocks in terms of
affine Hecke algebras (based on Bernstein's second adjointness
Theorem) have been given in many cases \cite{Hei}.
In this way affine Hecke algebras can be used for the explicit computation of
aspects of the harmonic analysis on $L$. The main result of this paper
will be another illustration of this point of view on affine Hecke algebras.

Analytic $R$-groups of affine Hecke algebras describe the
decomposition of tempered parabolically induced modules of the Hecke algebra
\cite{DeOp2}, in analogy with the theory of analytic $R$ groups for
reductive groups \cite{HC,KnSt,Sil}.
The analytic $R$-groups of affine Hecke algebras are defined in terms of
the Plancherel density, which is known known explicitly \cite{Opd-Sp,DeOp1,OpSo2}.
Therefore the analytic $R$-groups of affine Hecke algebras are amenable
to explicit determination \cite{Slo}. It is an interesting question how
the analytic $R$-groups for parabolic induction of tempered representations of
affine Hecke algebras are related to the analytic $R$-groups for parabolic
induction of tempered representations of reductive groups (some results in
this direction have been obtained by Roche \cite{Roc}). In the present paper
we will focus on the analytic $R$-groups of affine Hecke algebras,
in line with our point of view that affine Hecke algebras
are tools for explicit descriptions of the harmonic analysis on $L$.
\\[1mm]

For a $\C$-linear abelian category $\mathcal{C}$ with
finite homological dimension one defines the \emph{Euler-Poincar\'e pairing}
\cite{ScSt} of two objects of finite length
$\pi,\pi^\prime\in\textup{Obj}(\mathcal{C})$ by the formula
\begin{equation}\label{eq:ep}
EP(\pi,\pi^\prime)=
\sum_{i\geq 0}(-1)^i\textup{dim}\textup{Ext}_{\mathcal{C}}^i(\pi,\pi^\prime) .
\end{equation}
Let $K_\mathbb{C}(\mathcal{C})$ denote the Grothendieck group
(tensored by $\mathbb{C}$) of finite length objects in $\mathcal{C}$.
The Euler-Poincar\'e pairing extends to a sesquilinear form $EP$ on
$K_\C(\mathcal{C})$.

A case in point is the Euler-Poincar\'e pairing on the Grothendieck group
$G_\C (L)$ of admissible representations of $L$ \cite{BeDe,ScSt}.
In this situation the form $EP_L$ is Hermitian and plays a fundamental role in the
local trace formula and in the study of orbital integrals on the regular elliptic set of $L$
\cite{Art,ScSt,Bez,Ree}.
The definition of $EP$ also applies naturally to the Grothendieck group $G_\C (\mc{H})$
of finite dimensional representations of an affine Hecke algebra $\mc{H}$
with positive parameters. Here one uses that the category of finitely generated
$\mc{H}$-modules has finite cohomological dimension by \cite{OpSo1}.
The form $EP_\mc{H}$ on $G_\C (\mc{H})$ is Hermitian (\cite[Theorem 3.5 a)]{OpSo1}).

There is an obvious relation between these two instances of the form $EP$.
It follows first of all directly from the definition of $EP_L$ that the subspaces
$K_\C (\mathcal{B})$ of $G_\C (L)$ generated by the various
Bernstein blocks $\mathcal{B}$ of the category of smooth representations
of $L$ are mutually orthogonal with respect to the $EP$. Moreover, if
$\phi : \mc{B} \to \textup{Mod}(\mc{H})$ is an equivalence of categories then $\phi$
clearly induces a linear map from $K_\C (\mathcal{B})$ to
$G_\C (\mathcal{H})$ which respects the form $EP$.
The main result of the present paper is, for arbitrary abstract
affine Hecke algebras $\mc{H}$ with positive parameters,
the explicit computation of $EP_{\mc H}$ on the subspace of $G_\C (\mc{H})$
spanned by the finite dimensional tempered representations of $\mc{H}$,
in terms of the analytic $R$-groups of $\mc{H}$.

Our formula for $EP_\mc{H}$ is inspired by an analogous formula for $EP_L$ in terms of
the analytic $R$-groups of $L$, a result which emerges if one combines two deep results
in the harmonic analysis on $L$. The first ingredient is Kazhdan's orthogonality conjecture
\cite{Kaz} for $EP_L$, proved by \cite{ScSt} and \cite{Bez}. This result asserts that,
provided $\mathbb F$ has characteristic $0$ and the center of $L$ is compact, one has
\begin{equation}\label{eq:Kazhdan}
EP_L (\pi,\pi^\prime) =
\int_{C^{ell}}\theta_\pi (c^{-1})\theta_{\pi^\prime}(c) \textup{d} \gamma (c)
\end{equation}
for admissible representations $\pi,\pi^\prime$ of $L$. Here the space $C^{ell}$ denotes
the space of regular semisimple elliptic conjugacy classes of $L$, equipped with the canonical
``elliptic measure'' d$\gamma$ as defined by Kazhdan \cite{Kaz}, and $\theta_\pi$ denotes
the locally constant function on $C^{ell}$ determined by the distributional character of $\pi$.
Thus ``Kazhdan's conjecture'' expresses $EP(\pi,\pi^\prime)$ as an orthogonality pairing
of the characters of $\pi$ and $\pi^\prime$ with respect to the elliptic measure d$\gamma$
on the set $C^{ell}$.

The second ingredient is Arthur's formula expressing the right
hand side of \eqref{eq:Kazhdan} in terms of R-groups. More
precisely, let $P\subset L$ be a rational parabolic subgroup with
Levi component $M$, and let $\sigma$ be a smooth irreducible representation
of $M$, square integrable modulo the center of $M$. The representation
$\mathcal{I}_P^L(\sigma)$, the smooth normalized parabolically induced
representation from $\sigma$, is a tempered admissible unitary representation
of $L$. The decomposition of $\mathcal{I}_P^L(\sigma)$
is governed by the analytic R-group $\mf R_\sigma$, a finite group
which acts naturally on the real Lie algebra $\mr{Hom}(X^* (M), \R)$
of the center of $M$. For $r \in \mf R_\sigma$ we denote by $d(r)$
the determinant of the linear transformation $1 - r$ on
$\mr{Hom} (X^* (M), \R)$. Let $\pi$ be an irreducible tempered
representation of $L$ which occurs in $\mathcal{I}_P^L(\sigma)$
(denoted by $\pi\prec\mathcal{I}_P^L(\sigma)$).
The theory of the analytic R-group asserts that
$\rho=\textup{Hom}_L(\pi,\mathcal{I}_P^L(\sigma))$
is an irreducible projective $\mf R_\sigma$-representation.
We still assume that $\mathbb F$ has characteristic
zero and that the center of $L$ is compact. Applying the local
trace formula for $L$, Arthur \cite[Corollary 6.3]{Art} showed
that for all tempered irreducible representations
$\pi,\pi^\prime\prec \mathcal{I}_P^L(\sigma)$ one has
\begin{equation}\label{eq:Arthur}
\int_{C^{ell}}\theta_\pi (c^{-1})\theta_{\pi^\prime}(c) \textup{d} \gamma (c) =
|\mf R_\sigma|^{-1} \sum_{r \in \mf R_\sigma} |d(r)| \,
\overline{\mr{tr}_\rho (r)} \, \mr{tr}_{\rho'} (r)
\end{equation}
(In fact Arthur allowed $L$ to have noncompact center, but he had to
adjust \eqref{eq:Arthur} to make the integral converge.)
Combining \eqref{eq:Kazhdan} and \eqref{eq:Arthur}, one arrives at the
following formula \emph{provided that $\mathbb{F}$ has characteristic zero}:
\begin{equation}\label{eq:ArKaScSt}
EP_L (\pi,\pi') = |\mf R_\sigma|^{-1} \sum_{r \in \mf R_\sigma} |d(r)| \,
\overline{\mr{tr}_\rho (r)} \, \mr{tr}_{\rho'} (r) .
\end{equation}
This is the ``formula of Arthur'' referred to in the title of this paper. Before we assumed
that $L$ had compact center, but that is not necessary for \eqref{eq:ArKaScSt}.
If the center of $L$ is not compact, the equality is not very interesting however,
as both sides are zero. The left hand side by the argument of \cite[Lemma III.4.18.ii]{ScSt}
and the right hand side because every $r \in \mf R_\sigma$ fixes the subspace of
$\mr{Hom} (X^* (M),\R)$ corresponding to the $\mathbb F$-split part of the center of $L$
pointwise, so $d(r) = 0$.

The main results of the present paper are:

\begin{thmintro}\label{thm:1}(See Theorem \ref{thm:ArthurFormula}.) \\
The natural analog of equation (\ref{eq:ArKaScSt}) for $EP_\mc{H}$ for
affine Hecke algebras with arbitrary positive parameters holds true.
\end{thmintro}

In this result we use the notion of analytic R-groups for affine Hecke algebras
as introduced in \cite{DeOp2}.
The proof of Theorem \ref{thm:1} is completely different from the arguments sketched
above leading to equation (\ref{eq:ArKaScSt}), and is rather elementary in comparison to
the aforementioned deep results in $p$-adic harmonic analysis.
Although we wrote it down only in the case of affine Hecke algebras, this proof
of Theorem \ref{thm:1} applies to \emph{all} reductive $p$-adic groups.
In particular this observation proves:

\begin{thmintro}\label{thm:2} (See Remark \ref{rem:ArthurPadic}.) \\
Let $\mathbf{L}$ be a connected reductive group defined over an arbitrary
non-archimedian local field $\mathbb F$. Then \eqref{eq:ArKaScSt} holds for
$L = \mathbf{L}(\mathbb F)$.
\end{thmintro}

The elliptic integral in \eqref{eq:Kazhdan} and \eqref{eq:Arthur} does not seem to have an
obvious counterpart in the setting of affine Hecke algebras. Analogues are known for the
Euler--Poincar\'e pairing of representations of either a finite Weyl group \cite{Ree} or of an
(extended) affine Weyl group $W^e$ \cite[Theorem 3.3]{OpSo1}. To relate these to $EP_{\mc H}$
one must compare $EP_{\mc H}$ and $EP_{\mc W^e}$, which is done in \cite[Section 5.6]{Ree}
(for affine Hecke algebras with equal parameters) and in \cite[Chapter 3]{OpSo1}. Recent results
from \cite{Sol-Irr} allow us to conclude that $G_\C (\mc H)$ modulo the radical of $EP_{\mc H}$
equals the vector space Ell$(\mc H)$ of "elliptic characters", and that this space does not
depend on the parameters $q$ (Theorem \ref{thm:sigmaEll}).

Arthur's formula applies only to the Euler-Poincar\'e
pairing of tempered characters. In an abstract sense this no restriction, as it follows from
the Langlands classification of irreducible characters of $\mc{H}$ in terms of standard
induction data in Langlands position, that modulo the radical of the pairing $EP_{\mc H}$
any irreducible character is equivalent to a virtual tempered character. But this is
complicated in practice, and therefore our formula does not qualify as an explicit formula
for $EP_{\mc H}$ for general non-tempered irreducible characters.
It would therefore be desirable to extend the result to general finite dimensional
representations of $\mc{H}$. In the final section of this paper we make a first step by
extending the definition of the analytic $R$-group to non-tempered induction data.
We show that its irreducible characters are in natural bijection with the Langlands
quotients associated to the induction datum. However, we have not been able to
generalize the Arthur formula to this case.
\vspace{4mm}

\section{Preliminaries}

Here we recall the definitions and notations of our most important objects of study.
Several things described in this section can be found in more detail elsewhere in the
literature, see in particular \cite{Lus-Gr,Opd-Sp,OpSo1,Sol-Irr}.

\subsection{Affine Hecke algebras}

Let $\mf a$ be a finite dimensional real vector space and let $\mf a^*$ be its dual. Let
$Y \subset \mf a$ be a lattice and $X = \mr{Hom}_\Z (Y,\Z) \subset \mf a^*$ the dual lattice. Let
\[
\mc R = (X, R_0, Y ,R_0^\vee ,F_0) .
\]
be a based root datum. Thus $R_0$ is a reduced root system in $X ,\, R^\vee_0 \subset Y$
is the dual root system, $F_0$ is a basis of $R_0$ and the set of positive roots is denoted $R_0^+$.
Furthermore we are given a bijection $R_0 \to R_0^\vee ,\: \alpha \mapsto \alpha^\vee$ such
that $\inp{\alpha}{\alpha^\vee} = 2$ and such that the corresponding reflections
$s_\alpha : X \to X$ (resp. $s^\vee_\alpha : Y \to Y$) stabilize $R_0$ (resp. $R_0^\vee$).
We do not assume that $R_0$ spans $\mf a^*$.

The reflections $s_\alpha$ generate the Weyl group $W_0 = W (R_0)$ of $R_0$, and
$S_0 := \{ s_\alpha : \alpha \in F_0 \}$ is the collection of simple reflections. We have the
affine Weyl group $W^\af = \mh Z R_0 \rtimes W_0$ and the extended (affine) Weyl group
$W^e = X \rtimes W_0$. Both can be regarded as groups of affine transformations of
$\mf a^*$. We denote the translation corresponding to $x \in X$ by $t_x$.

As is well known, $W^\af$ is a Coxeter group, and the basis of $R_0$ gives rise to a set $S^\af$
of simple (affine) reflections. The length function $\ell$ of the Coxeter system $(W^\af ,S^\af )$
extends naturally to $W^e$. We write
\begin{align*}
&X^+ := \{ x \in X : \inp{x}{\alpha^\vee} \geq 0
\; \forall \alpha \in F_0 \} , \\
&X^- := \{ x \in X : \inp{x}{\alpha^\vee} \leq 0
\; \forall \alpha \in F_0 \} = -X^+ .
\end{align*}
It is easily seen that the center of $W^e$ is the lattice
\[
Z(W^e) = X^+ \cap X^- .
\]
We say that $\mc R$ is semisimple if $Z (W^e) = 0$ or equivalently if $R_0$ spans $\mf a^*$.
Thus a root datum is semisimple if and only if the corresponding reductive algebraic group is so.

With $\mc R$ we also associate some other root systems.
There is the non-reduced root system
\[
R_{nr} := R_0 \cup \{ 2 \alpha : \alpha^\vee \in 2 Y \} .
\]
Obviously we put $(2 \alpha )^\vee = \alpha^\vee / 2$. Let $R_1$
be the reduced root system of long roots in $R_{nr}$:
\[
R_1 := \{ \alpha \in R_{nr} : \alpha^\vee \not\in 2 Y \} .
\]
Consider a positive parameter function for $\mc R$, that is,
a function $q : W^e \to \R_{>0}$ such that
\begin{equation}\label{eq:parameterFunction}
\begin{array}{lll@{\quad}l}
q (\omega ) & = & 1 & \text{if } \ell (\omega ) = 0 , \\
q (w v) & = & q (w) q(v) & \text{if } w,v \in W^e \quad
\text{and} \quad \ell (wv) = \ell (w) + \ell (v) .
\end{array}
\end{equation}
Alternatively it can be given by $W_0$-invariant map $q : R_{nr}^\vee \to \R_{>0}$,
the relation being
\begin{equation}\label{eq:parameterEquivalence}
\begin{array}{lll}
q_{\alpha^\vee} = q(s_\alpha) = q (t_\alpha s_\alpha) & \text{if} & \alpha \in R_0 \cap R_1, \\
q_{\alpha^\vee} = q(t_\alpha s_\alpha) & \text{if} & \alpha \in R_0 \setminus R_1, \\
q_{\alpha^\vee / 2} = q(s_\alpha) q(t_\alpha s_\alpha)^{-1} & \text{if} &
\alpha \in R_0 \setminus R_1.
\end{array}
\end{equation}
The affine Hecke algebra $\mc H = \mc H (\mc R ,q)$ is the unique associative
complex algebra with basis $\{ N_w : w \in W^e \}$ and multiplication rules
\begin{equation}\label{eq:multrules}
\begin{array}{lll}
N_w \, N_v = N_{w v} & \mr{if} & \ell (w v) = \ell (w) + \ell (v) \,, \\
\big( N_s - q(s)^{1/2} \big) \big( N_s + q(s)^{-1/2} \big) = 0 & \mr{if} & s \in S^\af .
\end{array}
\end{equation}
In the literature one also finds this algebra defined in terms of the elements
$q(s)^{1/2} N_s$, in which case the multiplication can be described without square roots.
The algebra $\mc H$ is endowed with a conjugate-linear involution, defined on basis
elements by $N_w^* := N_{w^{-1}}$.

For $x \in X^+$ we put $\theta_x := N_{t_x}$. The corresponding semigroup morphism
$X^+ \to \mc H (\mc R ,q)^\times$ extends to a group homomorphism
\[
X \to \mc H (\mc R ,q)^\times : x \mapsto \theta_x .
\]
A part of the Bernstein presentation \cite[\S 3]{Lus-Gr} says that the subalgebra
$\mc A := \mr{span} \{ \theta_x : x \in X \}$ is isomorphic to $\C [X]$, and that the
center $Z (\mc H)$ corresponds to $\C [X]^{W_0}$ under this isomorphism.
Let $T$ be the complex algebraic torus
\[
T = \mr{Hom}_{\mh Z} (X, \mh C^\times ) \cong Y \otimes_\Z \C^\times ,
\]
so $\mc A \cong \mc O (T)$ and $Z (\mc H ) = \mc A^{W_0} \cong \mc O (T / W_0 )$.
This torus admits a polar decomposition
\[
T = T_{rs} \times T_{un} = \mr{Hom}_\Z (X, \R_{>0}) \times \mr{Hom}_\Z (X, S^1)
\]
into a real split (or positive) part and a unitary part.

Let $\mr{Mod}_f (\mc H)$ be the category of finite dimensional $\mc H$-modules,
and $\mr{Mod}_{f, W_0 t} (\mc H)$ the full subcategory of modules that admit the
$Z (\mc H)$-character $W_0 t \in T / W_0$. We let $G (\mc H)$ be the Grothendieck
group of $\mr{Mod}_f (\mc H)$ and we write $G_\C (\mc H) = \C \otimes_\Z G (\mc H)$.
Furthermore we denote by Irr$(\mc H)$, respectively $\mr{Irr}_{W_0 t} (\mc H)$,
the set of equivalence classes of irreducible objects in $\mr{Mod}_f (\mc H)$,
respectively $\mr{Mod}_{f, W_0 t} (\mc H)$.
We will use these notations also for other algebras and groups.

\subsection{Parabolic induction}

For a set of simple roots $P \subset F_0$ we introduce the notations
\begin{equation}\label{eq:parabolic}
\begin{array}{l@{\qquad}l}
R_P = \mh Q P \cap R_0 & R_P^\vee = \mh Q R_P^\vee \cap R_0^\vee , \\
X_P = X \big/ \big( X \cap (P^\vee )^\perp \big) &
X^P = X / (X \cap \mh Q P ) , \\
Y_P = Y \cap \mh Q P^\vee & Y^P = Y \cap P^\perp , \\
\mf a_P = \R P^\vee & \mf a^P = P^\perp , \\
T_P = \mr{Hom}_{\mh Z} (X_P, \mh C^\times ) &
T^P = \mr{Hom}_{\mh Z} (X^P, \mh C^\times ) , \\
\mc R_P = ( X_P ,R_P ,Y_P ,R_P^\vee ,P) & \mc R^P = (X,R_P ,Y,R_P^\vee ,P) .
\end{array}
\end{equation}
Although $T_{rs} = T_{P,rs} \times T_{rs}^P$, the product
$T_{un} = T_{P,un} T^P_{un}$ is not direct, because the intersection
\[
T_{P,un} \cap T_{un}^P = T_P \cap T^P
\]
can have more than one element (but only finitely many).

We define parameter functions $q_P$ and $q^P$ on the root
data $\mc R_P$ and $\mc R^P$, as follows. Restrict $q$ to a function on
$(R_P )_{nr}^\vee$ and use \eqref{eq:parameterEquivalence} to extend it to
$W^e (\mc R_P )$ and $W^e (\mc R^P )$. Now we can define the parabolic subalgebras
\[
\mc H_P = \mc H (\mc R_P ,q_P ) ,\qquad \mc H^P = \mc H (\mc R^P ,q^P ) .
\]
Despite our terminology $\mc H^P$ and $\mc H_P$ are not subalgebras of $\mc H$,
but they are close. Namely, $\mc H (\mc R^P ,q^P )$ is isomorphic to the subalgebra of
$\mc H (\mc R ,q)$ generated by $\mc A$ and $\mc H (W (R_P) ,q_P)$.
We denote the image of $x \in X$ in $X_P$ by $x_P$ and we let $\mc A_P \subset \mc H_P$
be the commutative subalgebra spanned by $\{ \theta_{x_P} : x_P \in X_P \}$.
There is natural surjective quotient map
\begin{equation}\label{eq:quotientP}
\mc H^P \to \mc H_P : \theta_x N_w \mapsto \theta_{x_P} N_w .
\end{equation}
For all $x \in X$ and $\alpha \in P$ we have
\[
x - s_\alpha (x) = \inp{x}{\alpha^\vee} \alpha \in \Z P,
\]
so $t (s_\alpha (x)) = t(x)$ for all $t \in T^P$. Hence $t (w(x)) = t (x)$ for all $w \in W (R_P)$,
and we can define an algebra automorphism
\begin{equation}
\phi_t : \mc H^P \to \mc H^P, \quad \phi_t (\theta_x N_w) = t (x) \theta_x N_w  \qquad t \in T^P .
\end{equation}
In particular, for $t \in T_P \cap T^P$ this descends to an algebra automorphism
\begin{equation}\label{eq:twistKP}
\psi_t : \mc H_P \to \mc H_P , \quad \theta_{x_P} N_w \mapsto t(x_P) \theta_{x_P} N_w
\qquad t \in T_P \cap T^P .
\end{equation}
Suppose that $g \in W_0$ satisfies $g (P) = Q \subseteq F_0$.
Then there are algebra isomorphisms
\begin{equation}\label{eq:psig}
\begin{array}{llcl}
\psi_g : \mc H_P \to \mc H_Q , &
\theta_{x_P} N_w & \mapsto & \theta_{g (x_P)} N_{g w g^{-1}} , \\
\psi_g : \mc H^P \to \mc H^Q , &
\theta_x N_w & \mapsto & \theta_{g x} N_{g w g^{-1}} .
\end{array}
\end{equation}

We can regard any representation $(\sigma ,V_\sigma)$ of $\mc H (\mc R_P ,q_P )$ as a
representation of $\mc H (\mc R^P ,q^P)$ via the quotient map \eqref{eq:quotientP}.
Thus we can construct the $\mc H$-representation
\[
\pi (P,\sigma ,t) := \mr{Ind}_{\mc H (\mc R^P ,q^P )}^{\mc H (\mc R ,q)} (\sigma \circ \phi_t ) .
\]
Representations of this form are said to be parabolically induced.

We intend to partition Irr$(\mc H )$ into finite packets, each of which is obtained
by inducing a discrete series representation of a parabolic subalgebra of $\mc H$.
The discrete series and tempered representations are defined via the
$\mc A$-weights of a representation, as we recall now.
Given $P \subseteq F_0$, we have the following positive cones in $\mf a$ and in $T_{rs}$:
\begin{equation}
\begin{array}{lll@{\qquad}lll}
\mf a^+ & = & \{ \mu \in \mf a : \inp{\alpha}{\mu} \geq 0 \:
  \forall \alpha \in F_0 \} , & T^+ & = & \exp (\mf a^+) , \\
\mf a_P^+ & = &  \{ \mu \in \mf a_P : \inp{\alpha}{\mu} \geq 0 \: \forall \alpha \in P \} , &
  T_P^+ & = & \exp (\mf a_P^+) , \\
\mf a^{P+} & = & \{ \mu \in \mf a^P : \inp{\alpha}{\mu} \geq 0 \:
  \forall \alpha \in F_0 \setminus P \} , & T^{P+} & = & \exp (\mf a^{P+}) , \\
\mf a^{P++} & = & \{ \mu \in \mf a^P : \inp{\alpha}{\mu} > 0 \:
  \forall \alpha \in F_0 \setminus P \} , & T^{P++} & = & \exp (\mf a^{P++}) .
\end{array}
\end{equation}
The antidual of $\mf a^{*+} :=  \{ x \in \mf a^* :  \inp{x}{\alpha^\vee} \geq 0
\: \forall \alpha \in F_0 \}$ is
\begin{equation}
\mf a^- = \{ \lambda \in \mf a : \inp{x}{\lambda} \leq 0 \: \forall x \in \mf a^{*+} \} =
\big\{ \sum\nolimits_{\alpha \in F_0} \lambda_\alpha \alpha^\vee : \lambda_\alpha \leq 0 \big\} .
\end{equation}
Similarly we define
\begin{equation}
\mf a_P^- = \big\{ \sum\nolimits_{\alpha \in P} \lambda_\alpha \alpha^\vee  \in \mf a_P :
\lambda_\alpha \leq 0 \big\} .
\end{equation}
The interior $\mf a^{--}$ of $\mf a^-$ equals
$\big\{ {\ts \sum_{\alpha \in F_0}} \lambda_\alpha \alpha^\vee : \lambda_\alpha < 0 \big\}$
if $F_0$ spans $\mf a^*$, and is empty otherwise. We write $T^- = \exp (\mf a^-)$ and
$T^{--} = \exp (\mf a^{--})$.

Let $t = |t| \cdot t |t|^{-1} \in T_{rs} \times T_{un}$ be the polar decomposition of $t$.
An $\mc H$-representation is called tempered if $|t| \in T^-$ for all its
$\mc A$-weights $t$, and anti-tempered if $|t|^{-1} \in T^-$ for all such $t$. More
restrictively we say that an irreducible  $\mc H$-representation belongs to the discrete
series (or simply: is discrete series) if $|t| \in T^{--}$, for all its $\mc A$-weights $t$.
In particular the discrete series is empty if $F_0$ does not span $\mf a^*$.

Our induction data are triples $(P,\delta,t)$, where
\begin{itemize}
\item $P \subset F_0$;
\item $(\delta ,V_\delta)$ is a discrete series representation of $\mc H_P$;
\item $t \in T^P$.
\end{itemize}
Let $\Xi$ be the space of such induction data, where we consider $\delta$ only modulo
equivalence of $\mc H_P$-representations. We say that $\xi = (P,\delta,t)$ is unitary
if $t \in T^P_{un}$, and we denote the space of unitary induction data by $\Xi_{un}$.
Similarly we say that $\xi$ is positive if $|t| \in T^{P+}$, which we write as $\xi \in \Xi^+$.
We have three collections of induction data:
\begin{equation}\label{eq:inductionData}
\Xi_{un} \subseteq \Xi^+ \subseteq \Xi .
\end{equation}
By default we endow these spaces with the topology for which $P$ and $\delta$ are
discrete variables and $T^P$ carries its natural analytic topology. With $\xi \in \Xi$
we associate the parabolically induced representation
\[
\pi (\xi) = \pi (P,\delta,t) := \mr{Ind}_{\mc H^P}^{\mc H} (\delta \circ \phi_t ) .
\]
As vector space underlying $\pi (\xi)$ we will always take $\C [W^P] \otimes V_\delta$,
where $W^P$ is the collection of shortest length representatives of $W_0 / W (R_P)$.
This space does not depend on $t$, which will allow us to speak of maps that
are continuous, smooth, polynomial or even rational in the parameter $t \in T^P$.
It is known that $\pi (\xi)$ is unitary and tempered if $\xi \in \Xi_{un}$, and non-tempered
if $\xi \in \Xi \setminus \Xi_{un}$, see \cite[Propositions 4.19 and 4.20]{Opd-Sp} and
\cite[Lemma 3.1.1]{Sol-Irr}.

The relations between such representations are governed by intertwining operators.
There construction \cite{Opd-Sp} is rather complicated, so we recall only their important
properties. Suppose that $P,Q \subset F_0 , k \in T_P \cap T^P , w \in W_0$ and $w (P) = Q$.
Let $\delta$ and $\sigma$ be discrete series representations of respectively $\mc H_P$ and
$\mc H_Q$, such that $\sigma$ is equivalent with $\delta \circ \psi_k^{-1} \circ \psi_w^{-1}$.

\begin{thm}\label{thm:intOp}
\textup{\cite[Theorem 4.33 and Corollary 4.34]{Opd-Sp}}
\enuma{
\item There exists a family of intertwining operators
\[
\pi (w k,P,\delta,t) : \pi (P,\delta,t) \to \pi (Q,\sigma,w (kt)) .
\]
As a map $\C [W^P] \otimes_\C V_\delta \to \C [W^Q] \otimes_\C V_\sigma$
it is a rational in $t \in T^P$ and constant on $T^{F_0}$-cosets.
\item This map is regular and invertible on an open neighborhood of $T^P_{un}$
in $T^P$ (with respect to the analytic topology).
\item $\pi (w k,P,\delta,t)$ is unitary if $t \in T^P_{un}$.
}
\end{thm}
We can gather all these intertwining operators in a groupoid $\mc W$, which we define now. The
base space of $\mc W$ is the power set of $F_0$ and the collection of arrows from $P$ to $Q$ is
\[
\mc W_{PQ} := \{ w \in W_0 : w (P) = Q \} \times T_P \cap T^P .
\]
Whenever it is defined, the multiplication in $\mc W$ is
\[
(w_1, k_1) \cdot (w_2, k_2) = (w_1 w_2, w_2^{-1}(k_1) k_2) .
\]

Families of intertwining operators $\pi (wk,P,\delta,t)$ satisfying the
properties listed in Theorem \ref{thm:intOp} are unique only up to normalization
by rational functions of $t\in T^P$ which are regular in an open neighborhood of
$T^P_{un}$, have absolute value equal to $1$ on $T^P_{un}$, and are constant on
$T^{F_0}$-cosets. The intertwining operators defined in \cite{Opd-Sp} are normalized
in such a way that composition of the intertwining operators corresponds to
multiplication of the corresponding elements of $\mc{W}$ only up to a scalar.

More precisely, let $g \in \mc W$ be such that $g w$ is defined. Then there exists a
number $\lambda \in \C$ with $|\lambda| = 1$, independent of $t$, such that
\begin{equation}\label{eq:multIntOp}
\pi (g ,Q,\sigma,w (k t)) \circ \pi (w k,P,\delta,t) =
\lambda \, \pi (g w k, P,\delta, t)
\end{equation}
as rational functions of $t$. We fix such a normalization of the intertwining
operators once and for all.

Let $W (R_P) r_\delta \in T_P / W (R_P)$ be the central character of the $\mc H_P$-representation
$\delta$. Then $|r_\delta| \in T_{P,rs} = \exp (\mf a_P)$, so we can define
\begin{equation}\label{eq:ccdelta}
cc_P (\delta) := W (R_P) \log |r_\delta| \in \mf a_P / W (R_P) .
\end{equation}
Since the inner product on $\mf a$ is $W_0$-invariant, the number $\norm{cc_P (\sigma)}$
is well-defined.

\begin{thm}\label{thm:parametrizationOfDual}
\textup{\cite[Theorem 3.3.2]{Sol-Irr}} \\
Let $\rho$ be an irreducible $\mc H$-representation. There exists a unique association class
$\mc W (P,\delta,t) \in \Xi / \mc W$ such that the following equivalent properties hold:
\enuma{
\item $\rho$ is isomorphic to an irreducible quotient of $\pi (\xi^+)$, for some
$\xi^+ = \xi^+ (\rho) \in \Xi^+ \cap \mc W (P,\delta,t)$;
\item $\rho$ is a constituent of $\pi (P,\delta,t)$, and $\norm{cc_P (\delta)}$
is maximal for this property.
}
\end{thm}

For any $\xi \in \Xi$ the packet
\begin{equation} \label{eq:packetReps}
\mr{Irr}_{\mc W \xi} (\mc H) := \{ \rho \in \mr{Irr} (\mc H) : \xi^+ (\rho) \in \mc W \xi \}
\end{equation}
is finite, but it is not so easy to predict how many equivalence classes of representations
it contains. This is one of the purposes of R-groups.

\subsection{The Schwartz algebra}

We recall how to complete an affine Hecke algebra to a topological algebra called the
Schwartz algebra. As a vector space $\mc S$ will consist of rapidly decreasing functions
on $W^e$, with respect to some length function. For this purpose it is unsatisfactory
that $\ell$ is 0 on the subgroup $Z(W^e)$, as this can be a large part
of $W$. To overcome this inconvenience, let
$L : X \otimes \mh R \to [0,\infty )$ be a function such that
\begin{itemize}
\item $L (X) \subset \mh Z$,
\item $L(x+y) = L(x) \quad \forall x \in X \otimes \mh R, y \in E$,
\item $L$ induces a norm on
$X \otimes \mh R / E \cong Z(W) \otimes \mh R$.
\end{itemize}
Now we define
\[
\mc N (w) := \ell (w) + L (w(0))  \qquad w \in W^e .
\]
Since $Z(W^e) \oplus \mh Z R_0$ is of finite index in $X$, the set
$\{ w \in W^e : \mc N (w) = 0 \}$ is finite. Moreover, because $W^e$ is
the semidirect product of a finite group and an abelian group, it is of
polynomial growth and different choices of $L$ lead to equivalent length
functions $\mc N$. For $n \in \mh N$ we define the norm
\[
p_n \big( {\ts \sum_{w \in W}} h_w N_w \big) :=
\sup_{w \in W} |h_w | (\mc N (w) + 1 )^n .
\]
The completion $\mc S = \mc S (\mc R ,q)$ of $\mc H (\mc R ,q)$ with respect to the family
of norms $\{ p_n \}_{n \in \mh N}$ is a nuclear Fr\'echet space. It consists of all possible
infinite sums $h = \sum_{w \in W} h_w N_w$ such that $p_n (h) < \infty \; \forall n \in \mh N$.
By \cite[Section 6.2]{Opd-Sp} or \cite[Appendix A]{OpSo2} $\mc S (\mc R,q)$ is a unital
locally convex *-algebra.

A crucial role in the harmonic analysis on affine Hecke algebra is played by a particular
Fourier transform, which is based on the induction data space $\Xi$. Let $\mc V_\Xi$
be the vector bundle over $\Xi$, whose fiber at $(P,\delta,t) \in \Xi$ is the representation
space $\C [W^P] \otimes V_\delta$ of $\pi (P,\delta,t)$. Let $\mr{End} (\mc V_\Xi)$ be the
algebra bundle with fibers $\mr{End}_\C (\C [W^P] \otimes V_\delta)$. Of course these
vector bundles are trivial on every connected component of $\Xi$, but globally not even the
dimensions need be constant. Since $\Xi$ has the structure of a complex algebraic variety,
we can construct the algebra of polynomial sections of $\mr{End} (\mc V_\Xi)$:
\[
\mc O \big( \Xi ; \mr{End} (\mc V_\Xi) \big) := \bigoplus_{P,\delta} \mc O (T^P) \otimes
\mr{End}_\C (\C [W^P] \otimes V_\delta) .
\]
For a submanifold $\Xi' \subset \Xi$ we define the algebra $C^\infty \big( \Xi' ;
\mr{End} (\mc V_\Xi) \big)$ in similar fashion.
The intertwining operators from Theorem \ref{thm:intOp} give rise to an action of the
groupoid $\mc W$ on the algebra of rational sections of $\mr{End} (\mc V_\Xi)$, by
\begin{equation}\label{eq:actSections}
(w \cdot f) (\xi) = \pi (w,w^{-1} \xi ) f (w^{-1} \xi) \pi (w, w^{-1} \xi )^{-1} ,
\end{equation}
whenever $w^{-1} \xi \in \Xi$ is defined. This formula also defines groupoid actions of
$\mc W$ on $C^\infty \big( \Xi' ; \mr{End} (\mc V_\Xi) \big)$, provided that $\Xi'$ is a
$\mc W$-stable submanifold of $\Xi$ on which all the intertwining operators are regular.
Given a suitable collection $\Sigma$ of sections of $(\Xi', \mr{End} (\mc V_\Xi) )$,
we write
\[
\Sigma^{\mc W} = \{ f \in \Sigma : (w \cdot f) (\xi) = f(\xi) \text{ for all } w \in \mc W,\
\xi \in \Xi' \text{ such that } w^{-1} \xi \text{ is defined} \} .
\]
The Fourier transform for $\mc H$ is the algebra homomorphism
\begin{align*}
& \mc F : \mc H \to \mc O \big( \Xi ; \mr{End} (\mc V_\Xi) \big) , \\
& \mc F (h) (\xi) = \pi (\xi) (h) .
\end{align*}
The very definition of intertwining operators shows that the image of $\mc F$ is contained
in the algebra $\mc O \big( \Xi ; \mr{End} (\mc V_\Xi) \big)^{\mc W}$.
By \cite[Theorem 5.3]{DeOp1} the Fourier transform extends to an isomorphism of
Fr\'echet *-algebras
\begin{equation}\label{eq:FourierIso}
\mc F : \mc S (\mc R,q) \to C^\infty \big( \Xi_{un} ; \mr{End} (\mc V_\Xi ) \big)^{\mc W} .
\end{equation}
Let $(P_1,\delta_1), \ldots, (P_N, \delta_N)$ be representatives for the action of $\mc W$ on
pairs $(P,\delta)$. Then the right hand side of \eqref{eq:FourierIso} can be rewritten as
\begin{equation}\label{eq:25}
\bigoplus\nolimits_{i=1}^N \big( C^\infty (T^P_{un}) \otimes
\mr{End} (\C [W^P] \otimes V_\delta) \big)^{\mc W_{P_i ,\delta_i}} ,
\end{equation}
where $\mc W_{P,\delta} = \{ w \in \mc W : w (P) = P, \delta \circ \psi_w^{-1} \cong \delta \}$
is the isotropy group of $(P,\delta)$.

\begin{rem}\label{rem:Plancherel}
Suppose that $L$ is a reductive $p$-adic group and that $K \subset L$ is a compact open
subgroup. Clearly \eqref{eq:FourierIso} is similar to the Plancherel isomorphism for the
Harish-Chandra--Schwartz algebra $\mc S (L)$ \cite{Wal}. One can easily deduce from \cite{Wal}
that the subalgebra $\mc S (L,K) \subset \mc S (L)$ of $K$-biinvariant functions has exactly the
same shape as \eqref{eq:25}, see \cite[Theorem 10]{Sol-Chern}.

Let $\mc H (L)$ be the Hecke algebra of $L$ and let Mod$(\mc H (L))$ be the category of smooth
$L$-representations. Consider the subalgebra $\mc H (L,K) = \mc S (L,K) \cap \mc H (L)$, for any
compact open subgroup $K \subset L$. According to \cite[Section 3]{BeDe} there exist arbitrarily
small $K$ such that Mod$(\mc H (L,K))$ is equivalent to the category consisting of those smooth
$L$-representations that are generated by their $K$-invariant vectors, and such that the latter
is a Serre subcategory of Mod$(\mc H (L))$. Moreover these subcategories exhaust Mod$(\mc H(L))$,
so all extensions of tempered smooth $L$-representations can be
described with representations of algebras $\mc S (L,K)$ of the form \eqref{eq:25}.
\end{rem}

It was shown in \cite[Corollary 5.5]{DeOp1} that \eqref{eq:FourierIso} implies
\begin{equation}
Z (\mc S) \cong C^\infty (\Xi_{un})^{\mc W},
\end{equation}
so the space of central characters of $\mc S$ is $\Xi_{un} / \mc W$.
We let $\mr{Mod}_f (\mc S)$ be the category of finite dimensional $\mc S$-modules
and $\mr{Mod}_{f,\mc W \xi} (\mc S)$ the full subcategory of modules which admit the
central character $\mc W \xi$. The collection $\mr{Irr}_{\mc W \xi} (\mc S)$ of
(equivalence classes of) irreducible objects in $\mr{Mod}_{f,\mc W \xi} (\mc S)$ equals
$\mr{Irr}_{\mc W \xi} (\mc H)$, in the notation of \eqref{eq:packetReps}.

For $\ep \in \mh R$ let $q^\ep$ be the parameter function $q^\ep (w) = q(w)^\ep$.
For every $\ep$ we have the affine Hecke algebra $\mc H (\mc R ,q^\ep )$ and its Schwartz
completion $\mc S (\mc R ,q^\ep )$. We note that $\mc H (\mc R ,q^0 ) = \mh C [W^e]$
is the group algebra of $W^e$ and that $\mc S (\mc R ,q^0 ) = \mc S (W^e)$ is the
Schwartz algebra of rapidly decreasing functions on $W^e$.

The intuitive idea is that these algebras depend continuously on $\ep$.
We will use this in the form of the following rather technical result.

\begin{thm}\label{thm:scalingReps}
\textup{\cite[Corollary 4.2.2]{Sol-Irr}} \\
For $\ep \in [-1,1]$ there exists a family of additive functors
\begin{align*}
& \tilde \sigma_\ep : \mr{Mod}_f (\mc H (\mc R ,q)) \to
\mr{Mod}_f (\mc H (\mc R ,q^\ep )) , \\
& \tilde \sigma_\ep (\pi ,V) = (\pi_\ep ,V) .
\end{align*}
with the properties
\begin{itemize}
\item[(1)] the map
\[
[-1,1] \to \mr{End}\, V : \ep \mapsto \pi_\ep (N_w)
\]
is analytic for any $w \in W^e$,
\item[(2)] $\tilde \sigma_\ep$ is a bijection if $\ep \neq 0$,
\item[(3)] $\tilde \sigma_\ep$ preserves unitarity,
\item[(4)] $\tilde \sigma_\ep$ preserves temperedness if $\ep \geq 0$,
\item[(5)] $\tilde \sigma_\ep$ preserves the discrete series if $\ep > 0$.
\end{itemize}
\end{thm}

\subsection{Analytic R-groups for tempered representations}
\label{sec:RgroupTemp}

We recall the definition of the analytic R-group from \cite{DeOp2}.
Let $\delta$ be a discrete series representation of $\mc H_P$ with central character
$W(R_P) r_\delta \in T_P / W(R_P)$. Recall that a \emph{parabolic subsystem}
of $R_0$ is a subset $R_Q \subset R_0$ satisfying $R_Q=R_0 \cap \mathbb{R}R_Q$.
We let $Q\subset R_{Q,+}:=R_Q\cap R_{0,+}$ denote the basis of $R_Q$ inside the
positive subset $R_{Q,+}\subset R_Q$. We call $R_Q$ \emph{standard} if
$Q\subset F_0$, in which case these notions agree with \eqref{eq:parabolic}.
The set of roots $R_0 \setminus R_P$ is a disjoint union
of subsets of the form $R_Q \setminus R_P$ where $R_Q\subset R_0$ runs over
the collection $\mc{P}^P_{min}$ of minimal parabolic subsystems
properly containing $R_P$. We define $\alpha_Q\in R_{0,+}$ by $Q=\{P,\alpha_Q\}$ and
$\alpha^P_Q=\alpha_Q \big|_{\mf a^P}\in \mf{a}^{P,*}$ for $R_Q\in\mc{P}^P_{min}$.
By the integrality properties of the root system $R_0$ we see that
\begin{equation}\label{eq:mult}
\{\alpha \big|_{\mf a^P}: \alpha\in R_Q\}\subset \mathbb{Z}\alpha_Q^P .
\end{equation}
Clearly $\alpha_Q^P$ is a character of $T^P$ whose kernel contains the
codimension one subtorus $T^Q\subset T^P$. For each $R_Q\in\mc{P}^P_{min}$
we define a $W (R_P)$-invariant rational function on $T$ by
\begin{equation}
c_Q^P(t):=\prod_{\alpha\in R_{Q,+} \setminus R_{P,+}}c_\alpha (t) ,
\end{equation}
where $c_\alpha$ denotes the usual rank one $c$-function associated to
$\alpha\in R_0$ and the parameter function $q$ of $\mc{H}$,
see e.g. \cite[Appendix 9]{DeOp1}.
For any $\alpha\in R_Q\backslash R_P$ the function $c_\alpha$ is a
nonconstant rational function on any coset of the form $r T^P$. In particular
$c_Q^P$ is regular on a nonempty Zariski-open subset of such a coset.
By the $W (R_P)$-invariance we see that for $t\in T^P$ and $r \in W (R_P ) r_\delta$,
the value of this function at $r t$ is independent of the choice
of $r \in W (R_P) r_\delta$. The resulting rational function $t \mapsto c_Q^P (r t)$
on $T^P$ is clearly constant along the cosets of
the codimension one subtorus $T^Q\subset T^P$.
We define the set of \emph{mirrors} $\mc{M}^{P,\delta}_Q$ associated
to $R_Q \in\mc{P}^P_{min}$ to be the set of connected components of the
intersection of the set of poles of this rational function with the unitary part
$T^P_{un}$ of $T^P$. We put
\begin{equation}
\mc{M}^{P,\delta}=\bigsqcup_{R_Q\in\mc{P}^P_{min}} \mc{M}^{P,\delta}_Q
\end{equation}
for the set of all mirrors in $T^P_{un}$ associated to the pair $(P,\delta)$.
Given $M\in\mc{M}^{P,\delta}$ we denote by $R_{Q^M}\subset R_0$ the
unique element of $\mc{P}^P_{min}$ such that $M\in \mc{M}_{Q^M}^{P,\delta}$.
Thus any mirror $M\in\mc{M}^{P,\delta}$ is a connected component of a
hypersurface of $T^P_{un}$ of the form $\alpha^P_{Q^M}=\textup{constant}$.
Observe that for a fixed pair $(P,\delta)$ the set $\mc{M}^{P,\delta}$ is finite
(and possibly empty).

For every mirror $M\in\mc{M}^{P,\delta}$ there exists,
by \cite[Theorem 4.3.i]{DeOp2}, a unique $s_M \in \mc{W}_{P,\delta}$ (i.e.
$s_M\in \mc{W}_{P,P}$ and $\delta\simeq\delta\circ\psi_{s_M}^{-1}$)
which acts as the reflection in $M$ on $T^P_{un}$.
For $\xi = (P,\delta,t) \in \Xi_{un}$ let $\mc M_\xi$ be the collection of mirrors
$M\in\mc{M}^{P,\delta}$ containing $t$. We define
\begin{equation}
R_\xi =  \{\pm \alpha^P_{Q^M} : M \in \mc M_\xi \}\text{ and }
R_\xi^+  =  \{\alpha^P_{Q^M} : M \in \mc M_\xi \} .
\end{equation}
Then it follows from \cite[Proposition 4.5]{DeOp2} and (\ref{eq:mult}) that $R_\xi$ is
a root system in $\mf{a}^{P,*}$, and $R_\xi^+\subset R_\xi$ is a positive subset.
Its Weyl group $W(R_{\xi})$ is generated by the reflections $s_M$ with $M \in \mc M_\xi$,
so it can be realized as a subgroup of $\mc W_\xi = \{ w \in \mc W : w (\xi) = \xi \}$.
The R-group of $\xi \in \Xi_{un}$ is defined as
\begin{equation}\label{eq:24}
\mf R_\xi = \{ w \in \mc W_\xi : w (R_{\xi}^+ ) = R_{\xi}^+ \} ,
\end{equation}
and by \cite[Proposition 4.7]{DeOp2} it is a complement to $W (R_{\xi})$ in $\mc W_\xi$:
\begin{equation}\label{eq:3}
\mc W_\xi = \mf R_\xi \ltimes W (R_{\xi}) .
\end{equation}
With these notions one can state the analogue of the Knapp--Stein linear independence theorem
\cite[Theorem 13.4]{KnSt} for affine Hecke algebras. For reductive $p$-adic groups the result
is proven in \cite{Sil}, see also \cite[Section 2]{Art}.

\begin{thm}\label{thm:KnappStein}
Let $\xi \in \Xi_{un}$.
\enuma{
\item For $w \in \mc W_\xi$ the intertwiner $\pi (w,\xi)$
is scalar if and only if $w \in W(R_{\xi})$.
\item There exists a 2-cocycle $\kappa_\xi$ (depending on the normalization of the
intertwining operators $\pi (w,\xi)$ for $w\in\mc{W}_{\xi}$) of $\mf R_\xi$
such that $\mr{End}_{\mc H} (\pi (\xi))$
is isomorphic to the twisted group algebra $\C [\mf R_\xi ,\kappa_\xi]$.
\item Given the normalization of the intertwining operators, there is a unique bijection
\[
\mr{Irr}(\C [\mf R_\xi ,\kappa_\xi]) \to \mr{Irr}_{\mc W \xi} ( \mc S ) :
\rho \mapsto \pi_\rho ,
\]
such that $\pi (\xi) \cong \bigoplus_\rho \pi_\rho \otimes \rho$ as
$\mc H \otimes \C [\mf R_\xi ,\kappa_\xi]$-modules.
\item Let $\mr{Mod}_{f,un,\mc W \xi}(\mc S)$ be the category of finite dimensional unitary
$\mc S$-representations that admit the central character $\mc W \xi$. The functor
\[
\begin{array}{lccc}
E_\xi : & \mr{Mod}_f (\C [\mf R_\xi ,\kappa_\xi]) & \to & \mr{Mod}_{f, un, \mc W \xi} (\mc S) ,\\
& \rho & \mapsto & \mr{Hom}_{\C [\mf R_\xi ,\kappa_\xi]} (\rho ,\pi (\xi)) .
\end{array}
\]
is an equivalence of categories, with inverse $E_\xi^{-1}(\pi) = \mr{Hom}_{\mc H}(\pi, \pi (\xi))$.
}
\end{thm}
\noindent \emph{Proof.}
Part (a) is \cite[Theorem 5.4]{DeOp2}, parts (b) and (c) are \cite[Theorem 5.5]{DeOp2} and
part (d) is \cite[Theorem 5.13]{Opd-ICM}. $\qquad \Box$
\vspace{4mm}

\section{The Euler-Poincar\'e pairing}

Let $\pi, \pi' \in \mr{Mod}_f (\mc H)$. Their Euler-Poincar\'e pairing is defined as
\begin{equation}\label{eq:4}
EP_{\mc H} (\pi,\pi') = \sum_{n \geq 0} (-1)^n \dim_\C \mr{Ext}_{\mc H}^n (\pi,\pi') .
\end{equation}
Since $\mc H$ is Noetherian and has finite cohomological dimension \cite[Proposition 2.4]{OpSo1},
this pairing is well-defined.

\begin{thm}\label{thm:EP}
\enuma{
\item For all $\pi,\pi' \in \mr{Mod}_f (\mc H)$ and $\ep \in [0,1] $:
\[
EP_{\mc H (\mc R, q^\ep)} (\tilde \sigma_\ep (\pi), \tilde \sigma_\ep (\pi')) = EP_{\mc H} (\pi,\pi') .
\]
\item The pairing $EP_{\mc H}$ is symmetric and positive semidefinite.
\item If $P \subset F_0$ and $\R P \neq \mf a^*$, then $\mr{Ind}_{\mc H^P}^{\mc H}(V)$
lies in the radical of $EP_{\mc H}$, for all $V \in \mr{Mod}_f (\mc H^P)$.
}
\end{thm}
\noindent \emph{Proof.}
See Proposition 3.4 and Theorem 3.5 of \cite{OpSo1}. We note that the symmetry
of $EP_{\mc H}$ is not automatic, it is proved via part (a) for $\ep = 0$ and a more detailed
study of $EP_{W^e} = EP_{\mc H (\mc R ,q^0)}$ \cite[Theorem 3.2]{OpSo1}.
$\qquad \Box$ \\[2mm]

The pairing $EP_{\mc H}$ extends naturally to a Hermitian form on $G_\C (\mc H)$, say complex
linear in the second argument. By Theorem \ref{thm:EP}.c it factors through the quotient
\[
\mr{Ell} (\mc H) := G_\C (\mc H) \Big/ \sum_{P \subset F_0, \R P \neq \mf a^*}
\mr{Ind}_{\mc H^P}^{\mc H} \big( G_\C (\mc H^P) \big)
\]
Following \cite{Opd-ICM} we call Ell$(\mc H)$ the space of elliptic characters of $\mc H$. Notice
that Ell$(\mc H) = 0$ and $EP_{\mc H} = 0$ if the root datum $\mc R$ is not semisimple.

\begin{lem}\label{lem:Elltemp}
The composite map $G_\C (\mc S) \to G_\C (\mc H) \to \mr{Ell} (\mc H)$ is surjective.
\end{lem}
\noindent \emph{Proof.} We may and will assume that $\mc R$ is semisimple.
We have to show that every irreducible $\mc H$-representation $\pi$ that does not vanish
in Ell$(\mc H)$, can be written as a linear combination of tempered representations and of
representations induced from proper parabolic subalgebras.
Recall \cite[Section 2.2]{Sol-Irr} that a Langlands datum is a triple $(P,\sigma,t)$ such that
\begin{itemize}
\item $P \subset F_0$ and $\sigma$ is an irreducible tempered $\mc H_P$-representation;
\item $t \in T^P$ and $|t| \in T^{P++}$.
\end{itemize}
The Langlands classification \cite[Theorem 2.2.4]{Sol-Irr} says that the $\mc H$-representation
$\pi (P,\sigma,t)$ has a unique irreducible quotient $L(P,\sigma,t)$ and that there is
(up to equivalence) a unique Langlands datum such that $L(P,\sigma,t) \cong \pi$.

If $P = F_0$, then $T^P = \{1\}$ because $\mc R$ is semisimple.
So $\mc H_P = \mc H^P = \mc H$ and $\pi \cong \sigma$, which by definition is tempered.

Therefore we may suppose that $P \neq F_0$. Since $\pi (P,\sigma,t)$ is induced from $\mc H^P$,
we have $\pi = L(P,\sigma,t) - \pi (P,\sigma,t)$ in Ell$(\mc H)$. By \cite[Lemma 2.2.6.b]{Sol-Irr}
\[
\pi (P,\sigma,t) - L (P,\sigma,t) \in G_\C (\mc H)
\]
is a sum of representations $L(Q,\tau,s)$ with $P \subset Q$ and
$\norm{cc_P (\sigma)} < \norm{cc_Q (\tau)}$, where $cc_P (\sigma)$ is as in
\eqref{eq:ccdelta}. Given the central character of $\pi$ (an element of $T/W_0$),
there are only finitely many possibilities for $cc_P (\sigma)$. Hence we can deal with
the representations $L(Q,\tau,s)$ via an inductive argument. $\qquad \Box$
\\[2mm]

From \cite[Lemma 4.2.3.a]{Sol-Irr} we know that $\tilde \sigma_0 : \mr{Mod}_f (\mc H) \to
\mr{Mod}_f (W^e)$ commutes with parabolic induction, so it induces a linear map
\begin{equation}
\sigma_{\mr{Ell}} : \mr{Ell}(\mc H) \to \mr{Ell}(W^e) = \mr{Ell} \big( \mc H (\mc R ,q^0) \big) .
\end{equation}
\begin{thm}\label{thm:sigmaEll}
The pairings $EP_{\mc H}$ and $EP_{W^e}$ induce Hermitian inner products on respectively
Ell$(\mc H)$ and Ell$(W^e)$, and the map $\sigma_{\mr{Ell}}$ is an isometric bijection.
\end{thm}
\noindent \emph{Proof.}
By Lemma \ref{lem:Elltemp} and \cite[(3.35)]{Sol-Irr}, $\sigma_{\mr{Ell}}$ is a linear bijection.
According to \cite[Theorem 3.2.b]{OpSo1} $EP_{W^e}$ induces a Hermitian inner product on
Ell$(W^e)$, and by Theorem \ref{thm:EP} $\sigma_{\mr{Ell}}$ is an isometry. Therefore the
sesquilinear form on Ell$(\mc H)$ induced by $EP_{\mc H}$ is also a Hermitian inner product.
$\qquad \Box$
\vspace{4mm}
\section{The Arthur formula}
\label{sec:Arthur}

In this section we prove the Arthur formula for the Euler--Poincar\'e pairing \eqref{eq:4}
of representations of affine Hecke algebras. In view of Lemma \ref{lem:Elltemp} little is lost
if we restrict $EP_{\mc H}$ to tempered representations, so we do that until further notice.
One of the main results of \cite{OpSo1} says that
\begin{equation}\label{eq:ExtHS}
\mr{Ext}^n_{\mc S} (\pi,\pi') \cong \mr{Ext}^n_{\mc H} (\pi,\pi') \qquad n \in \Z_{\geq 0}
\end{equation}
for all $\pi,\pi' \in \mr{Mod}_f (\mc S)$, and hence
\begin{equation} \label{eq:5}
EP_{\mc S} (\pi,\pi') := \sum_{n \geq 0} (-1)^n \dim_\C \mr{Ext}_{\mc S}^n (\pi,\pi')
\quad \text{equals} \quad EP_{\mc H}(\pi,\pi') .
\end{equation}
We remark that corresponding statement for reductive $p$-adic groups
is a deep theorem of Meyer \cite{Mey}.

It is clear from the definition \eqref{eq:4} that $EP_{\mc H}(\pi,\pi') = 0$ if $\pi$ and $\pi'$
admit different $Z(\mc H)$-characters. With \eqref{eq:5} we can strengthen this to
$EP_{\mc H}(\pi,\pi') = 0$ whenever $\pi$ and $\pi'$ admit different $Z(\mc S)$-characters.
This is really stronger, because $Z(\mc S)$ is larger than the closure of $Z(\mc H)$ in $\mc S$.

We fix $\xi = (P,\delta,t) \in \Xi_{un}$. As $\Xi_{un}$ has the structure of a smooth
manifold (with components of different dimensions) we can consider its tangent space
$T_\xi (\Xi_{un})$ at $\xi$. The isotropy group $\mc W_\xi$ acts linearly on $T_\xi (\Xi_{un})$
and for $w \in \mc W_\xi$ we denote the determinant of the corresponding linear map by
$\det (w )_{T_\xi (\Xi_{un})} $. Since the connected component of $\xi$ in $\Xi_{un}$ is
diffeomorphic to $T^P_{un}$, we have
\[
T_\xi (\Xi_{un}) \cong \text{Lie}(T^P_{un}) \cong i \mf a^P \cong i (\mf a / \mf a_P) .
\]
The expression $\det (1-w )_{T_\xi (\Xi_{un})} = \det (1-w )_{\mf a^P}$
is analogous to the Weyl denominator and to $d(r)$ in Theorem \ref{thm:1}. Notice that
\[
\det (1-w )_{T_\xi (\Xi_{un})} \geq 0
\]
because $T_\xi (\Xi_{un})$ is a real representation of the finite group $\mc W_\xi$. Clearly
$\det (1-w )_{T_\xi (\Xi_{un})} \neq 0$ if and only if $w$ acts without fixed points on
$T_\xi (\Xi_{un}) \setminus \{0\}$, in which case we say that $w$ is elliptic in $\mc W_\xi$.

Given $\pi \in \mr{Mod}_{f,un,\mc W \xi} (\mc S)$, Theorem \ref{thm:KnappStein}.d
produces a representation $E_\xi^{-1} (\pi) = \mr{Hom}_{\mc H}(\pi ,\pi (\xi))$ of
$\C [\mf R_\xi, \kappa_\xi]$. Explicitly, this means
\begin{equation}\label{eq:1}
E_\xi^{-1} (\pi) (w_1) E_\xi^{-1} (\pi) (w_2) E_\xi^{-1} (\pi) (w_1 w_2 )^{-1} =
\kappa_\xi (w_1,w_2) \quad \text{for all } w_1,w_2 \in \mc W_\xi ,
\end{equation}
A variation on \eqref{eq:1} shows that the dual space $E_\xi^{-1} (\pi)^*$ is a representation
of $\C [\mf R_\xi, \kappa_\xi^{-1}]$. If $\pi' \in \mr{Mod}_{f, un, \mc W \xi} (\mc S)$, then
$E_\xi^{-1} (\pi)^* \otimes E_\xi^{-1} (\pi')$ is a representation of $\C [\mf R_\xi,1]$, that is,
of the group $\mf R_\xi$. Its trace is
\begin{equation}\label{eq:23}
\mr{tr}_{E_\xi^{-1} (\pi)^* \otimes E_\xi^{-1} (\pi')} (r) =
\overline{\mr{tr}_{E_\xi^{-1} (\pi)} (r)} \, \mr{tr}_{E_\xi^{-1} (\pi')} (r) .
\end{equation}
It is an elementary result in homological algebra that, for the purpose of computing
Euler--Poincar\'e-pairings, one may replace any module by its semisimplification.
Hence in the next theorem it suffices to compute $EP_\mc{H}(\pi,\pi')$ for
$\pi,\pi'\in \mr{Mod}_{f,\mc W \xi}(\mc S)$ irreducible. Recall that irreducible
tempered modules are unitarizable \cite[Corollary 3.23]{DeOp1}. In particular Theorem
\ref{thm:KnappStein}.d applies to $\pi$ and $\pi'$.
\begin{thm}\label{thm:ArthurFormula}
Let $\pi, \pi' \in \mr{Mod}_{f,\mc W \xi} (\mc S) $ be irreducible tempered $\mc{H}$-modules.
Denote by $\rho:=E_\xi^{-1}(\pi)$ and $\rho':=E_\xi^{-1}(\pi')$
the irreducible modules of $\mathbb{C}[\mf{R}_\xi,\kappa_\xi]$ corresponding to
$\pi$ and $\pi'$ respectively (according to Theorem \ref{thm:KnappStein}).
Then we have:
\[
EP_{\mc H} (\pi,\pi') = | \mf R_\xi |^{-1} \sum_{r \in \mf R_\xi} \det (1-r )_{T_\xi
(\Xi_{un})}\overline{\mr{tr}_{\rho} (r)} \, \mr{tr}_{\rho'} (r) .
\]
\end{thm}
\noindent \emph{Proof.}
Our strategy is to consider $\pi$ and $\pi'$ as modules over various algebras $A$, such that
the extension groups $\mr{Ext}_A^n (\pi, \pi')$ for different $A$'s are all isomorphic.

Write $\xi = (P,\delta,t)$. By \eqref{eq:5}, \eqref{eq:FourierIso} and \eqref{eq:25} we may
replace $\mc H$ first by $\mc S$ and then by
\begin{equation}\label{eq:2}
\big( C^\infty (T^P_{un}) \otimes \mr{End}_\C (\C [W^P] \otimes V_\delta ) \big)^{\mc W_{P,\delta}} .
\end{equation}
Let $U_\xi \subset T^P_{un}$ be a nonempty open ball around $t$ such that
\[
\overline{U_\xi} \cap w \overline{U_\xi} = \left\{ \begin{array}{lll}
\overline{U_\xi} & \text{if} & w \in \mc W_\xi \\
\emptyset & \text{if} & w \in \mc W_{P,\delta} \setminus \mc W_\xi .
\end{array} \right.
\]
The localization of \eqref{eq:2} at $U_\xi$ is
\begin{equation}
\begin{split}\label{eq:6}
C^\infty (U_\xi)^{\mc W_\xi} \otimes_{C^\infty (T^P_{un})^{\mc W_{P,\delta}}}
\big( C^\infty (T^P_{un}) \otimes \mr{End}_\C (\C [W^P] \otimes V_\delta ) \big)^{\mc W_{P,\delta}} \\
\cong \big( C^\infty (U_\xi) \otimes \mr{End}_\C (\C [W^P] \otimes V_\delta ) \big)^{\mc W_\xi} .
\end{split}
\end{equation}
As localization is an exact functor, this preserves the Ext-groups, up to a natural isomorphism.
The action of $\mc W_\xi$ on $C^\infty (U_\xi) \otimes \mr{End}_\C (\C [W^P] \otimes V_\delta )$
is still defined by \eqref{eq:actSections}, so for $w \in \mc W_\xi$ and $t' \in U_\xi$
\begin{equation}\label{eq:7}
(w \cdot f)(t') = \pi (w,P,\delta, w^{-1}t') f (w^{-1} t') \pi (w,P,\delta, w^{-1}t')^{-1} .
\end{equation}
Let $V_\xi$ be the vector space $\C [W^P] \otimes V_\delta$ endowed with the
$\mc H$-representation $\pi (\xi)$. It is also a projective $\mc W_\xi$-representation, so we can
define a $\mc W_\xi$-action on $C^\infty (U_\xi) \otimes \mr{End}_\C (V_\xi)$ by
\begin{equation}\label{eq:8}
(w \cdot f)(t') = \pi (w,\xi) f (w^{-1} t') \pi (w,\xi)^{-1} .
\end{equation}
The difference between \eqref{eq:7} and \eqref{eq:8} is that $\pi (w,P,\delta, w^{-1}t')$
depends on $t'$, while $\pi (w,\xi)$ does not. Since $U_\xi$ is $\mc W_\xi$-equivariantly
contractible to $t$, we are in the right position to apply \cite[Lemma 7]{Sol-Chern}. Its proof, and
in particular \cite[(20)]{Sol-Chern} shows that the algebra \eqref{eq:6}, with the $\mc W_\xi$-action
\eqref{eq:7}, is isomorphic to
\begin{equation}\label{eq:9}
\big( C^\infty (U_\xi) \otimes \mr{End}_\C (V_\xi ) \big)^{\mc W_\xi},
\end{equation}
with respect to the $W_\xi$-action \eqref{eq:8}.
By Theorem \ref{thm:KnappStein}.a the elements of $W (R_{\xi})$ do not really act on
$\mr{End}_\C (V_\xi)$, only on $C^\infty (U_\xi)$. Moreover by \eqref{eq:3} $W (R_{\xi})$ is
normal in $\mc W_\xi$ and $W_\xi / W (R_{\xi}) \cong \mf R_\xi$, so we can rewrite \eqref{eq:9} as
\begin{equation}\label{eq:10}
\big( C^\infty (U_\xi )^{W (R_{\xi})} \otimes \mr{End}_\C (V_\xi ) \big)^{\mf R_\xi} .
\end{equation}
Using suitable coordinates on $U_\xi$, we can identify it as a $\mc W_\xi$-representation with
$T_\xi (\Xi_{un})$. By Chevalley's Theorem \cite{Che,Hum} the algebra
$\mc O (T_\xi (\Xi_{un}))^{W (R_{\xi})}$ is a free polynomial algebra on $\dim T^P_{un}$
generators. With \cite[Theorem 3]{Bie} or \cite{Schw}
one can extend Chevalley's Theorem to smooth functions. Explicitly, this means that
$U_\xi / W(R_{\xi})$ is a manifold with corners such that
\begin{equation}\label{eq:11}
C^\infty (U_\xi )^{W (R_{\xi})} \cong C^\infty \big( U_\xi / W(R_{\xi}) \big)
\cong C^\infty \big( \R^{\dim T^P_{un}} \big) .
\end{equation}
We insert \eqref{eq:11} in \eqref{eq:10} to obtain
\begin{equation}\label{eq:12}
\big( C^\infty \big( U_\xi / W(R_{\xi}) \big) \otimes \mr{End}_\C (V_\xi ) \big)^{\mf R_\xi} .
\end{equation}
Let $\kappa_\xi : \mf R_\xi \times \mf R_\xi \to \C^\times$ be the cocycle from Theorem
\ref{thm:KnappStein}.b and let $\mf R_\xi^*$ be a Schur extension of $\mf R_\xi$
(also known as a representation group, see \cite[Section 53]{CuRe}).
There is a central idempotent $p \in \C [\mf R_\xi^*]$ such that
\[
\C [\mf R_\xi ,\kappa_\xi] \cong p \C [\mf R_\xi^*]
\]
as algebras and as $\C [\mf R_\xi]$-bimodules. By Theorem \ref{thm:KnappStein}.c every
irreducible representation of $\C [\mf R_\xi,\kappa_\xi]$ appears at least
once in $V_\xi$, so by Theorem \ref{thm:Moritaeq}.a \eqref{eq:12} is Morita equivalent to
\begin{equation}\label{eq:16}
\big( C^\infty \big( U_\xi / W(R_{\xi}) \big) \otimes
\mr{End}_\C (p \C [\mf R^*_\xi] ) \big)^{\mf R^*_\xi} .
\end{equation}
Under this Morita equivalence $\pi$ corresponds to a representation
$\pi_1$ of \eqref{eq:16} such that
\begin{equation}\label{eq:17}
\begin{split}
\rho := \mr{Hom}_{\mc H} (\pi, \pi (\xi)) \cong \mr{Hom}_{\big( C^\infty \big( U_\xi /
W(R_{\xi}) \big) \otimes \mr{End}_\C (V_\xi ) \big)^{\mf R_\xi}} (\pi, \pi (\xi)) \\
\cong \mr{Hom}_{\big( C^\infty \big( U_\xi / W(R_{\xi}) \big) \otimes
\mr{End}_\C (p \C [\mf R^*_\xi] ) \big)^{\mf R^*_\xi}} (\pi_1, p \C [\mf R^*_\xi]) ,
\end{split}
\end{equation}
as $\mf R^*_\xi$-representations. Since $p$ is central and by Theorem \ref{thm:Moritaeq}.b,
\eqref{eq:16} is isomorphic to
\begin{equation}\label{eq:18}
p \big( C^\infty \big( U_\xi / W(R_{\xi}) \big) \otimes \mr{End}_\C (\C [\mf R^*_\xi] )
\big)^{\mf R^*_\xi} \cong p \big( C^\infty (U_\xi / W(R_{\xi})) \rtimes \mf R^*_\xi \big) .
\end{equation}
This is a direct summand of $C^\infty (U_\xi / W(R_{\xi})) \rtimes \mf R^*_\xi$, so we can
consider the Ext-groups of $\pi_1$ and $\pi'_1$ just as well with respect to
$C^\infty (U_\xi / W(R_{\xi})) \rtimes \mf R^*_\xi$.
More precisely, by Theorem \ref{thm:Moritaeq}.c and \eqref{eq:17} we get
\begin{equation}\label{eq:21}
\mr{Ext}_{\mc H}^n (\pi,\pi') \cong \mr{Ext}^n_{C^\infty (U_\xi / W(R_{\xi}))
\rtimes \mf R^*_\xi} (\pi_2,\pi'_2) ,
\end{equation}
where $\pi_2$ admits the $C^\infty (U_\xi / W(R_{\xi}))$-character $\mf R^*_\xi t$
and $\pi_2 \big|_{\mf R^*_\xi} \cong \rho^*$, and similarly for $\pi'_2$ and $\rho'$.

Now that we took enough reduction steps, we can resort to more direct calculations. The simple
structure of $C^\infty (U_\xi / W(R_{\xi})) \rtimes \mf R^*_\xi$ allows one to construct an
analogue of the Koszul resolution, which was used in \cite[Theorem 3.2.c]{OpSo1} to show that
\begin{equation}\label{eq:19}
EP_{C^\infty (U_\xi / W(R_{\xi})) \rtimes \mf R^*_\xi} (\pi_2,\pi'_2) =
e_{\mf R^*_\xi} (\pi_2, \pi'_2) = e_{\mf R^*_\xi} (\rho^*, \rho'^*) .
\end{equation}
Following \cite{Ree}, $e_{\mf R^*_\xi}$ denotes the elliptic pairing of
$\mf R^*_\xi$-representations, with respect to the representation of $\mf R^*_\xi$ on the
tangent space $T_t (U_\xi / W(R_{\xi}))$ of $U_\xi / W(R_{\xi})$ at $t$. Of course this is imprecise,
as the orbifold $U_\xi / W(R_{\xi})$ has a corner at $t$. To make more sense of it, one uses
\eqref{eq:11} and reinterprets $T_t (U_\xi / W(R_{\xi}))$ as the tangent space of
$\R^{\dim T^P_{un}}$ at 0. These identifications can be made $\mf R_\xi$-equivariantly,
so $T_0 \big( \R^{\dim T^P_{un}} \big)$ becomes a $\mf R_\xi$-representation.

According to \cite[Section 2]{Ree}, the right hand side of \eqref{eq:19} equals
\begin{equation}\label{eq:20}
\begin{split}
| \mf R^*_\xi |^{-1} \sum\nolimits_{r \in \mf R^*_\xi} \det (1-r )_{T_t (U_\xi / W(R_{\xi}))}
\overline{\mr{tr}_{\pi_2}(r)} \mr{tr}_{\pi'_2}(r) \; = \\
| \mf R^*_\xi |^{-1} \sum\nolimits_{r \in \mf R^*_\xi} \det (1-r )_{T_t (U_\xi / W(R_{\xi}))}
\mr{tr}_{\rho}(r) \overline{\mr{tr}_{\rho'}(r)} .
\end{split}
\end{equation}
Notice that this formula does not use the entire action of $\mf R^*_\xi$ on
$T_t (U_\xi / W(R_{\xi}))$, only $\det (1-r )_{T_t (U_\xi / W(R_{\xi}))}$. This determinant is zero
whenever $r \in \mf R^*_\xi$ fixes a nonzero vector in $T_t (U_\xi / W(R_{\xi}))$. Suppose that
$R_{\xi}$ is nonempty. Then $R_{\xi}^\vee$ is a nonempty root system in $\mf a^P$ and
$\R R_{\xi}^\vee$ can be identified with a subspace of $T_t (U_\xi)$. Every $r \in \mf R^*_\xi$
fixes $\sum_{\alpha \in R^+_{\xi}} \alpha^\vee \in \mf a^P \setminus \{0\}$, so $r$ fixes
nonzero vectors of $T_t (U_\xi)$ and $T_t (U_\xi / W(R_{\xi}))$. We conclude that
\[
\det (1-r )_{T_t (U_\xi / W(R_{\xi}))} = 0 = \det (1-r )_{T_t (U_\xi)}
\]
whenever $R_{\xi}$ is nonempty. Therefore we may always replace
$T_t (U_\xi / W(R_{\xi}))$ by $T_t (U_\xi) = T_\xi (\Xi_{un})$ in \eqref{eq:20}.
So by \eqref{eq:21}, \eqref{eq:19} and \eqref{eq:20}
\begin{equation}\label{eq:22}
EP_{\mc H} (\pi, \pi') = | \mf R^*_\xi |^{-1} \sum_{r \in \mf R^*_\xi}
\det (1-r )_{T_\xi (\Xi_{un})} \mr{tr}_{\rho}(r) \overline{\mr{tr}_{\rho'}(r)} .
\end{equation}
Finally we want to reduce from $\mf R^*_\xi$ to $\mf R_\xi$. The action of $\mf R^*_\xi$ on
$T_\xi (\Xi_{un})$ is defined via the quotient map $\mf R^*_\xi \to \mf R_\xi$, so that is no problem.
We noted in \eqref{eq:23} that $\rho\otimes\rho'^*$ is a $\mf R_\xi$-representation with trace
\[
\mr{tr}_{\rho\otimes\rho'^*} (r) =
\mr{tr}_\rho (r) \overline{\mr{tr}_{\rho'}(r)} .
\]
Hence any two elements of $\mf R^*_\xi$ with the same image in $\mf R_\xi$ give the same
contribution to \eqref{eq:22}. Furthermore $EP_{\mc H}$ is symmetric by Theorem \ref{thm:EP}.b, so
\begin{equation*}
EP_{\mc H} (\pi, \pi') = EP_{\mc H} (\pi' ,\pi) = | \mf R_\xi |^{-1} \sum_{r \in \mf R_\xi}
\det (1-r )_{T_\xi (\Xi_{un})} \mr{tr}_{\rho'}(r) \overline{\mr{tr}_\rho (r)} .
\quad \Box
\end{equation*}

\begin{rem}\label{rem:ArthurPadic}
The proof of Theorem \ref{thm:ArthurFormula} relies on three deep results: the Knapp--Stein
linear independence theorem (Theorem \ref{thm:KnappStein}), the isomorphism
$\mr{Ext}^*_{\mc H} \cong \mr{Ext}^*_{\mc S}$ for tempered representations \eqref{eq:ExtHS}
and the Plancherel isomorphism for $\mc S$ \eqref{eq:FourierIso}. Since these three theorems are
also known for all reductive $p$-adic groups (as we indicated in some remarks) the proof of Theorem
\ref{thm:ArthurFormula} applies equally well to the Euler--Poincar\'e pairing for admissible
tempered representations of reductive $p$-adic groups. This proves Theorem \ref{thm:2}.
\end{rem}
\vspace{4mm}

\section{Algebras of invariants}

Let $U$ be a smooth manifold, possibly with corners, and let $G$ be a finite group
acting on $U$ by diffeomorphisms $\alpha_g$. Let $V$ be a finite dimensional complex
$G$-representation and endow the algebra $C^\infty (U) \otimes \mr{End}_\C (V)$
with the diagonal $G$-action. Algebras of the form
\[
\big( C^\infty (U) \otimes \mr{End}_\C (V) \big)^G
\]
were studied among others in \cite{Sol-Chern} and \cite[Section 2.5]{Sol-Thesis}.
In the proof of Theorem \ref{thm:ArthurFormula} we used some elementary results
on such algebras, which we prove here. Certainly a substantial part of this section is
already known, but the authors have not found a suitable reference.

\begin{thm}\label{thm:Moritaeq}
Let $V'$ be another finite dimensional $G$-representation and assume that, for every irreducible
$G$-representation $W, \; \mr{Hom}_G (W,V) = 0$ if and only if $\mr{Hom}_G (W,V') = 0$.
\enuma{
\item The algebras $\big( C^\infty (U) \otimes \mr{End}_\C (V) \big)^G$ and
$\big( C^\infty (U) \otimes \mr{End}_\C (V') \big)^G$ are Morita-equivalent.
\item For the regular representation $\C [G]$ of $G$ we have
\[
\big( C^\infty (U) \otimes \mr{End}_\C (\C [G]) \big)^G \cong C^\infty (U) \rtimes G .
\]
\item Denote the isotropy group of $u \in U$ by $G_u$. Let $\sigma$ be finite dimensional
representation $G_u$ and let $\sigma^*$ be the contragredient representation.
Under the isomorphism (b) the $C^\infty (U) \rtimes G$-representation
\[
\mr{Ind}_{C^\infty (U) \rtimes G_u}^{C^\infty (U) \rtimes G} (\C_u \otimes \sigma)
\]
corresponds to a $\big( C^\infty (U) \otimes \mr{End}_\C (\C [G]) \big)^G$-representation
$\pi (u,\sigma)$ with $C^\infty (U)^G$-character $G u$, such that
\[
\mr{Hom}_{\big( C^\infty (U) \otimes \mr{End}_\C (V) \big)^G}
(\pi (u,\sigma), \C_u \otimes V)) \cong \sigma^*,
\]
as $G_u$-representations.
}
\end{thm}
\noindent \emph{Proof.}
(a) Consider the bimodules
\[
\begin{array}{lll}
M_1 & = & C^\infty (U) \otimes \mr{Hom}_\C (V,V'), \\
M_2 & = & C^\infty (U) \otimes \mr{Hom}_\C (V',V).
\end{array}
\]
It is clear that
\begin{equation}\label{eq:13}
M_1 \otimes_{C^\infty (U) \otimes \mr{End}_\C (V)} M_2 \cong C^\infty (U) \otimes \mr{End}_\C (V') ,
\end{equation}
and similarly in the reversed order. We will show that \eqref{eq:13} remains valid if we take
$G$-invariants everywhere. In other words, we claim that the bimodules $M_1^G$ and $M_2^G$
implement the desired Morita-equivalence. Notice the multiplication yields a natural map
\begin{equation}\label{eq:14}
M_1^G  \otimes_{\big( C^\infty (U) \otimes \mr{End}_\C (V) \big)^G} M_2^G  \to
\big( C^\infty (U) \otimes \mr{End}_\C (V') \big)^G .
\end{equation}
Take any $u \in U$ and let $I \subset C^\infty (U)^G$ be the maximal ideal of functions
vanishing at $G u$. Dividing out \eqref{eq:14} by $I$, we obtain
\begin{equation}\label{eq:15}
\mr{Hom}_{G_u} (V,V') \otimes_{\mr{End}_{G_u} (V')} \mr{Hom}_{G_u} (V',V) \to
\mr{End}_{G_u} (V) .
\end{equation}
Decompose the $G_u$-representations $V$ and $V'$ as
\[
V \cong \bigoplus_{\rho \in \mr{Irr}(G)} \rho \otimes \C^{m_\rho} ,\qquad
V' \cong \bigoplus_{\rho \in \mr{Irr}(G)} \rho \otimes \C^{m'_\rho}
\]
The assumption of the Theorem says that $m_\rho = 0$ if and only if $m'_\rho = 0$. As
\[
\mr{End}_{G_u} (V') \cong
\bigoplus_{\rho \in \mr{Irr}(G)} \mr{End}_{G_u} (\rho) \otimes \mr{End}_\C (\C^{m'_\rho}) ,
\]
the left hand side of \eqref{eq:15} is isomorphic to
\begin{multline*}
\hspace{-5mm} \bigoplus_{\rho \in \mr{Irr}(G)} \!\! \mr{End}_{G_u} (\rho)
\otimes \mr{Hom}_\C (\C^{m_\rho}, \C^{m'_\rho}) \: \otimes_{\mr{End}_{G_u} (V')}
\bigoplus_{\rho \in \mr{Irr}(G)} \mr{End}_{G_u} (\rho) \otimes
\mr{Hom}_\C (\C^{m'_\rho}, \C^{m_\rho}) \\
\cong \bigoplus_{\rho \in \mr{Irr}(G)} \!\! \mr{End}_{G_u} (\rho) \otimes
\mr{End}_\C (\C^{m_\rho}) \cong \mr{End}_{G_u} (V) .
\end{multline*}
Hence \eqref{eq:15} is a bijection for all $u \in U$, which implies that \eqref{eq:14} is injective.
The image of \eqref{eq:14} is a two-sided ideal of $\big( C^\infty (U) \otimes \mr{End}_\C (V')
\big)^G$, which is dense in the sense that for every $u \in U$ the algebra and its ideal have the
same set of values at $u$. Consequently \eqref{eq:14} is bijective, and it is an isomorphism of
$\big( C^\infty (U) \otimes \mr{End}_\C (V') \big)^G$-bimodules.\\
(b) For use in part (c) we provide a proof of this folklore result.
We agree that the action of $G$ on $V = \C [G]$ is $\rho(g) (h) = h g^{-1}$.
For a decomposable tensor $f \otimes g \in C^\infty (U) \rtimes G$ we define
$L(f \otimes g) \in C^\infty (U) \otimes \mr{End}_\C (\C [G])$ by
\[
L(f \otimes g) (h) = \alpha_{h^{-1} g^{-1}} (f) \otimes gh \qquad h \in G .
\]
It is easily verified that $L(f \otimes g)$ is $G$-invariant and that $L$
extends to an algebra homomorphism
\[
L : C^\infty (U) \rtimes G \to \big( C^\infty (U) \otimes \mr{End}_\C (\C [G]) \big)^G .
\]
We claim that $L$ is invertible, with inverse
\[
L'(b) = \sum_{g \in G} b(g^{-1})_e \otimes g \text{, where }
b(g^{-1}) = \sum_{h \in G} b (g^{-1})_h \otimes h \in C^\infty (U) \otimes \C [G] .
\]
Clearly $L' (L (f \otimes g)) = f \otimes g$, while for any
$b \in \big( B \otimes \mr{End}(\mh C [G]) \big)^G$ and $h \in G$:
\begin{equation}
\begin{split}
L(L'b)(h) &= \sum_{g \in G} \alpha_{h^{-1} g^{-1}} (b(g^{-1})_e) \otimes gh \\
 &= \alpha_{h^{-1}} \big( (\alpha_{g^{-1}} b)(g^{-1}) \big)_e \otimes gh \\
 &= \sum_{g \in G} \alpha_{h^{-1}} \big( (\rho(g) b \rho (g^{-1}) ) (g^{-1}) \big)_e \otimes gh \\
 &= \sum_{g \in G} \alpha_{h^{-1}} \big( b(e) g^{-1} \big)_e \otimes gh \\
 &= \big( \alpha_{h^{-1}} b \big) (e) h \\
 &= \big( \rho(h^{-1}) \alpha_{h^{-1}}(b) \rho(h) \big) (h) \;=\; b(h) .
\end{split}
\end{equation}
Thus indeed $L'$ is the inverse of $L$.\\
(c) It suffices to show this when $\sigma$ is an irreducible $G_u$-representation.
The isomorphism from part (b) preserves the central characters, so the ideal $I$ of
$\big( C^\infty (U) \otimes \mr{End}_\C (V) \big)^G$ generated by
$\{ f \in C^\infty (U)^G | f(u) = 0 \}$ annihilates $\pi (u,\sigma)$. In other words, we
must describe $\pi (u,\sigma)$ as a representation of
\[
\big( C^\infty (U) \otimes \mr{End}_\C (V) \big)^G / I \cong
\big( \bigoplus\nolimits_{g \in G / G_u} \mr{End}_\C (\C [G]) \big)^G \cong \mr{End}_{G_u} (\C [G]) .
\]
As a $G_u$-representation $V = \C [G] \cong \C [G_u] \otimes \C^{[G:G_u]}$, so
\[
\mr{End}_{G_u} (\C [G]) \cong \mr{End}_{G_u} (\C [G_u]) \otimes \mr{End}_\C (\C^{[G:G_u]})
\cong \C [G_u] \otimes M_{[G:G_u]} (\C) .
\]
Unter these identifications $\pi (u,\sigma)$ corresponds to the representation
$\sigma \otimes \C^{[G:G_u]}$. Then
\[
\mr{Hom}_{\big( C^\infty (U) \otimes \mr{End}_\C (V) \big)^G} (\pi (u,\sigma), \C_u \otimes V)
\cong \mr{Hom}_{G_u} (\sigma, \C [G_u]) .
\]
As a $G_u \times G_u$-representation $\C [G_u] \cong \bigoplus_{\rho \in \mr{Irr}(G_u)}
\rho \otimes \rho^*$, so $\mr{Hom}_{G_u} (\sigma, \C [G_u]) \cong \sigma^*$
as $G_u$-representations. $\qquad \Box$
\vspace{4mm}

\section{Analytic R-groups for non-tempered representations}

In this section we extend the definitions and results from Section \ref{sec:RgroupTemp}
(which stem from \cite{DeOp2}) to non-tempered induction data $\xi \in \Xi \setminus \Xi_{un}$.
One approach would be to define the R-group of $\xi = (P,\delta,t)$ to be that of
$\xi_{un} = (P,\delta ,t \,|t|^{-1})$. This works out well in most cases, but it is unsatisfactory
if the isotropy group $\mc W_{\xi_{un}}$ is strictly larger than $\mc W_\xi$.

The problem is that $t \, |t|^{-1}$ is not always a generic point of $(T^P)^{\mc W_\xi}$.
To overcome this, we define $\mc M_\xi$ to be the collection of mirrors $M \in \mc M^{P,\delta}$
that contain not only $t \, |t|^{-1}$, but the entire connected component of $t |t|^{-1}$ in
$(T^P_{un})^{\mc W_\xi}$. Notice that all components of $(T^P)^{\mc W_\xi}$ are cosets of a
complex subtorus of $T^P$. Since $t \, |t|^r \in (T^P)^{\mc W_\xi}$ for all $r \in \R ,\; \mc M_\xi$
is precisely the collection of mirrors whose complexification contains $t$.

Now we define $R_{\xi} ,R^+_{\xi}$ and $\mf R_\xi$ as in \eqref{eq:24} and \eqref{eq:3}.
We note that $\mf R_\xi$ may differ from $\mf R_{\xi_{un}}$, but that $R_\xi = R_{P,\delta,t'}$
for almost all $t' \in (T^P_{un})^{\mc W_\xi}$. For positive induction data there is an analogue
of Theorem \ref{thm:KnappStein}:

\begin{thm} \label{thm:KnappSteinNontemp}
Let $\xi = (P,\delta,t) \in \Xi^+$.
\enuma{
\item For $w \in \mc W_\xi$ the intertwiner $\pi (w,\xi)$ is scalar if and only if $w \in W(R_{\xi})$.
\item There exists a 2-cocycle $\kappa_\xi$ of $\mf R_\xi$ such that
$\mr{End}_{\mc H} (\pi (\xi))$ is isomorphic to the twisted group algebra $\C [\mf R_\xi ,\kappa_\xi]$.
\item There exists a unique bijection
\[
\mr{Irr}(\C [\mf R_\xi ,\kappa_\xi]) \to \mr{Irr}_{\mc W \xi} ( \mc H ) :
\rho \mapsto \pi_\rho ,
\]
and there exist indecomposable direct summands $\widetilde{\pi_\rho}$ of $\pi (\xi)$,
such that $\pi_\rho$ is a quotient of $\widetilde{\pi_\rho}$ and
$\pi (\xi) \cong \bigoplus_\rho \widetilde{\pi_\rho} \otimes \rho$ as
$\mc H \otimes \C [\mf R_\xi ,\kappa_\xi]$-modules.
}
\end{thm}
\begin{rem} Part (a) also holds for general $\xi \in \Xi$. Part (d) of Theorem
\ref{thm:KnappStein} does not admit a nice generalization to $\xi \in \Xi^+$.
One problem is that the category $\mr{Mod}_f (\C [\mf R_\xi ,\kappa_\xi])$ is semisimple,
while $\pi (\xi)$ is not always completely reducible. One can try to consider the category of
$\mc H$-representations all whose irreducible subquotients lie in $\mr{Irr}_\xi (\mc H)$,
but then one does not get a nice formula for the functor $E_\xi$.
\end{rem}
\emph{Proof.}
(a) and (b) Since $\xi$ is positive, the operators $\pi (w,\xi)$ with $w \in \mc W_\xi$ span
$\mr{End}_{\mc H} (\pi (\xi))$ \cite[Theorem 3.3.1]{Sol-Irr}.
Let $v \in W(R_{\xi})$. The intertwiner $\pi (v,P,\delta ,t')$ is rational in $t'$ and
by Theorem \ref{thm:KnappStein}.a it is scalar for all $t'$ in a Zariski-dense subset of
$(T^P_{un})^{\mc W_\xi}$. Hence $\pi (v,P,\delta ,t')$ is scalar for all
$t' \in (T^P)^{\mc W_\xi}$, which together with \eqref{eq:3} and \eqref{eq:multIntOp} implies that
\[
\mr{End}_{\mc H} (\pi (\xi)) = \mr{span} \{ \pi (w,\xi) : w \in \mf R_\xi \} .
\]
All the intertwiners $\pi (w,P,\delta,t')$ depend continuously on $t'$, so the type of
$\pi (P,\delta,t')$ as a projective $\mf R_\xi$-representation is constant on connected
components of $(T^P)^{\mc W_\xi}$. Again by Theorem \ref{thm:KnappStein}.a
$\{ \pi (w,\xi) : w \in \mf R_\xi \}$ is linearly independent for generic
$t' \in (T^P_{un})^{\mc W_\xi}$, so it is in fact linearly independent for all $t' \in (T^P)^{\mc W_\xi}$.
Now the multiplication rules for intertwining operators \eqref{eq:multIntOp} show that
$\mr{End}_{\mc H} (\pi (\xi))$ is isomorphic to a twisted group algebra of $\mf R_\xi$.

(c) Let $A = \mr{End}_{\C [\mf R_\xi,\kappa_\xi]}$ be the bicommutant of $\pi (\xi, \mc H)$ in
$\mr{End}_\C (\mc H)$. There is a canonical bijection
\[
\mr{Irr}(\C [\mf R_\xi ,\kappa_\xi]) \to \mr{Irr} (A) : \rho \mapsto \widetilde{\pi_\rho}
\]
such that $\pi (\xi) \cong \bigoplus_\rho \widetilde{\pi_\rho} \otimes \rho$ as
$A \otimes \C [\mf R_\xi ,\kappa_\xi]$-modules. By construction the irreducible $A$-subrepresentations
of $\pi (\xi)$ are precisely the indecomposable $\mc H$-subrepresentations of $\pi (\xi)$.
By \cite[Proposition 3.1.4]{Sol-Irr} every $\widetilde{\pi_\rho}$ has a unique irreducible quotient
$\mc H$-representation, say $\pi_\rho$, and $\pi_\rho \cong \pi_{\rho'}$ if and only if
$\widetilde{\pi_\rho} \cong \widetilde{\pi_{\rho'}}$ as $\mc H$-representations. Thus
$\rho \mapsto \pi_\rho$ sets up a bijection between Irr$(\C [\mf R_\xi ,\kappa_\xi])$ and the
equivalence classes of irreducible quotients of $\pi (\xi)$. By Theorem \ref{thm:parametrizationOfDual}
the latter collection is none other than $\mr{Irr}_{\mc W \xi} (\mc H). \qquad \Box$

\end{document}